\documentclass[11pt]{amsart} 
\usepackage{amssymb,amsmath,latexsym,enumerate,graphicx,bbm,mathptmx,microtype}
\usepackage{pgf, tikz}
\usetikzlibrary{decorations.pathreplacing}
\usetikzlibrary{matrix,arrows} 

\newtheorem{theorem}{Theorem} 
\newtheorem{corollary}[theorem]{Corollary}

\newtheorem{exam}{Example}

\newtheorem*{rem}{Remarks}

\def\mult{\operatorname{mult}}
\def\Z{\mathbb{Z}}

\def\R{\mathbb{R}}

\def\HH{\mathcal{H}}
\def\P{\mathcal{P}}

\def\m{\mathbf{m}}

\def\x{\mathbf{x}}

\begin{document}

\title{A bivariate chromatic polynomial for signed graphs}

\author{Matthias Beck}
\address{Department of Mathematics\\
         San Francisco State University\\
         1600 Holloway Ave\\
         San Francisco, CA 94132\\
         U.S.A.}
\email{mattbeck@sfsu.edu}

\author{Mela Hardin}
\address{School of Mathematical and Statistical Sciences\\
		Arizona State University\\
		Physical Sciences Building A-Wing Rm. 216\\
		901 S. Palm Walk\\
		Tempe, AZ 85287\\
		U.S.A.}
\email{melahardin@asu.edu}

\begin{abstract}
We study Dohmen--P\"onitz--Tittmann's bivariate chromatic polynomial $c_\Gamma(k,l)$ which counts all $(k+l)$-colorings of a graph $\Gamma$ such that adjacent vertices get different colors if they are $\le k$.
Our first contribution is an extension of $c_\Gamma(k,l)$ to signed graphs, for which we obtain an inclusion--exclusion formula and several special evaluations giving rise, e.g., to polynomials that encode balanced subgraphs.
Our second goal is to derive combinatorial reciprocity theorems for $c_\Gamma(k,l)$ and its signed-graph analogues, reminiscent of Stanley's reciprocity theorem linking chromatic polynomials to acyclic orientations.
\end{abstract}



\date{7 June 2014}

\thanks{We thank the referees for helpful suggestions. This research was partially supported by the U.\ S.\ National Science Foundation through the grants DMS-1162638 (Beck) and DGE-0841164 (Hardin).}

\maketitle


\section{Introduction}

Graph coloring problems are ubiquitous in many areas within and outside of mathematics.
For a positive integer $n$, let $[n] := \left\{ 1, 2, \dots, n \right\}$ and $[\pm n] := \left\{ -n, -n+1, \dots, n \right\}$.
We study the \emph{bivariate chromatic polynomial} $c_\Gamma(k,l)$ of a graph $\Gamma = (V, E)$, first introduced in \cite{dohmenponitztittmann} and defined as the counting function of colorings $\x \in [k+l]^V$ that satisfy for any edge $vw \in E$
\[
  x_v \ne x_w \qquad \text{ or } \qquad x_v = x_w > k \, .
\]
The usual chromatic polynomial of $\Gamma$ can be recovered as the special evaluation $c_\Gamma(k,0)$.
Dohmen, P\"onitz, and Tittmann provided basic properties of $c_\Gamma(k,l)$ in \cite{dohmenponitztittmann}, including polynomiality and special evaluations yielding the matching and independence polynomials of $\Gamma$.
Subsequent results include a deletion--contraction formula and applications to Fibonacci-sequence identities \cite{hillarwindfeldt}, common generalizations of $c_\Gamma(k,l)$ and the Tutte polynomial \cite{averbouchgodlinmakowsky}, and closed formulas for paths and cycles \cite{dohmenbivariatepathsandcycles}.

Our first goal is to introduce and study the natural analogue of the bivariate chromatic polynomial for signed graphs, which originated in the social sciences and have found applications also in biology, physics, computer science, and economics; see \cite{zaslavskydynamicsurvey} for a comprehensive bibliography.
A \emph{signed graph} $\Sigma = (\Gamma,\sigma)$ consists of a graph $\Gamma = (V,E)$ and a signature $\sigma \in \left\{ \pm \right\}^E$.
The underlying graph $\Gamma$ may have multiple edges and, besides the usual links and loops, also \emph{halfedges} (with only one endpoint) and \emph{loose edges} (no endpoints); the latter are irrelevant for coloring questions, and so we assume in this paper that $\Sigma$ has no loose edges.
An unsigned graph can be realized by a signed graph all of whose edges are labelled with $+$.

We define the function $c_\Sigma(2k+1,2l)$ as counting the \emph{proper $(k,l)$-colorings} $\x \in [\pm (k+l)]^V$, namely, those colorings that satisfy for any edge $vw \in E$
\[
  x_v \ne \sigma_{vw} \, x_w \qquad \text{ or } \qquad |x_v| = |x_w| > k \, .
\]
This \emph{bivariate chromatic polynomial} (in Corollary \ref{chromaticpolcor} we will see that $c_\Sigma(2k+1,2l)$ is indeed a polynomial) specializes to Zaslavsky's chromatic polynomial of signed graphs \cite{zaslavskysignedcoloring} in the case $l=0$.
As in Zaslavsky's theory, $c_\Sigma(2k+1,2l)$ comes with a companion, the \emph{zero-free bivariate chromatic polynomial} $c^*_\Sigma(2k,2l)$ which counts all proper $(k,l)$-colorings $\x \in ([\pm (k+l)] \setminus 0)^V$.

Our first result is a deletion--contraction formula, the common generalization of \cite[Theorem 2.3]{zaslavskysignedcoloring} and \cite[Lemma 1.1]{hillarwindfeldt}. The definitions of deletions and contractions of signed graphs are reviewed in detail in Section~\ref{bivariateresultssection}, where we also prove our other results for the bivariate chromatic polynomials.

\begin{theorem}\label{Spolynomial}
Let $\Sigma$ be a signed graph.
If $e\in E$ is a halfedge or negative loop then
\[
  c_\Sigma(2k+1,2l) = c_{\Sigma - e} (2k+1,2l) - c_{\Sigma/e} (2k+1,2l) \, ;
\]
if $e\in E$ is not a halfedge or negative loop then
\[
  c_\Sigma(2k+1,2l) = c_{\Sigma - e} (2k+1,2l) - c_{\Sigma/e} (2k+1,2l) + 2l \, c_{[\Sigma/e] - v} (2k+1,2l)
\]
and
\[
  c^*_\Sigma(2k,2l) = c^*_{\Sigma - e} (2k,2l) - c^*_{\Sigma/e} (2k,2l) + 2l \, c^*_{[\Sigma/e] - v} (2k,2l) \, ,
\]
where $v$ is the vertex to which $e$ contracts in $\Sigma/e$.
\end{theorem}

A component of the signed graph $\Sigma = (\Gamma, \sigma)$ is \emph{balanced} if it contains no halfedges and each cycle has positive sign product, and it is \emph{antibalanced} if its negative $(\Gamma, -\sigma)$ is balanced.
We define the \emph{antibalance polynomial} of $\Sigma$ as
\[
  a_\Sigma (x,y) := \sum_{ S \text{ antibalanced subgraph of } \Sigma } x^{ |V(S)| } y^{ c(S) } ,
\]
where $c(S)$ denotes the number of components of $S$.
This polynomial relates to the zero-free bivariate chromatic polynomial as follows:

\begin{theorem}\label{chromaticantibalancethm}
$c^*_\Sigma (2, 2l) = (2l)^{ |V| } \, a_\Sigma ( \frac{ 1 }{ 2l } , 2) \, .$
\end{theorem}

Our second goal is to prove \emph{reciprocity theorems} for the bichromatic polynomials for graphs and signed graphs, in analogy with the following theorem of Stanley \cite{stanleyacyclic} on the usual chromatic polynomial $c_\Gamma (k)$.
We call an orientation of a graph $\Gamma$ \emph{compatible} with the coloring $\x \in \Z^V$ if $x_v \le x_w$ for any edge oriented from $v$ to~$w$.

\begin{theorem}[Stanley]\label{stanleyrecthm}
For $k \in \Z_{ >0 }$,
$(-1)^{ |V| } c_\Gamma (-k)$ equals the number of $k$-colorings of $\Gamma$, each counted with multiplicity equal to the number of compatible acyclic orientations of $\Gamma$.
In particular, $(-1)^{ |V| } c_\Gamma (-1)$ equals the number of acyclic orientations of~$\Gamma$.
\end{theorem}

Our generalization for bivariate chromatic polynomials is as follows.

\begin{theorem}\label{bivariaterecthm}
For $k \in \Z_{ >0 }$ and $l \in \Z_{ \ge 0 }$,
$(-1)^{ |V| } c_\Gamma (-k,-l)$ equals the number of $(k+l)$-colorings of $\Gamma$, each counted with multiplicity: a $k$-coloring has multiplicity equal to the number of compatible acyclic orientations of $\Gamma$, and a coloring that uses at least one color $> k$ has multiplicity~$1$.
\end{theorem}

We prove this theorem in Section \ref{reciprocitysection}, where we also give an analogous reciprocity theorem for the bivariate chromatic polynomials of signed graphs.
We finish with a few open problems in Section~\ref{sec:openprob}.


\section{Bivariate Chromatic Polynomials for Signed Graphs}\label{bivariateresultssection}

We first review a few constructs on a signed graph $\Sigma = (V, E, \sigma)$.
The \emph{restriction} of $\Sigma$ to an edge set $F \subseteq E$ is the signed graph $\left( V, F, \sigma|_F \right)$.
For $e \in E$, we denote by $\Sigma - e$ (the \emph{deletion} of $e$) the restriction of $\Sigma$ to $E - \{ e \}$.
For $v \in V$, denote by $\Sigma - v$ the restriction of $\Sigma$ to $E-F$ where $F$ is the set of all edges incident to~$v$.

\emph{Switching} $\Sigma$ by $s \in \left\{ \pm \right\}^V$ results in the new signed graph $(V, E, \sigma^s)$ where
$
  \sigma^s_{ vw } = s_v \, \sigma_{ vw } \, s_w
$.
Switching does not alter balance, and any balanced signed graph can be obtained from switching an all-positive graph \cite{zaslavskysignedgraphs}.
We also note that there is a natural bijection of colorings of $\Sigma$ and a switched version of it, and this bijection preserves the number of proper $(k,l)$-colorings.

The \emph{contraction} of $\Sigma$ by $F \subseteq E$, denoted by $\Sigma/F$, is defined as follows \cite{zaslavskysignedgraphs}:
switch $\Sigma$ so that every balanced component of $F$ is all positive, coalesce all vertices of each balanced component, and discard the remaining vertices and all edges in $F$; note that this may produce halfedges.
For example, if $F = \left\{ e \right\}$ for a link $e$, $\Sigma/e$ is obtained by switching $\Sigma$ so that $\sigma(e) = +$ and then contracting $e$ as in the case of unsigned graphs, that is, disregard $e$ and identify its two endpoints.
As a second example, if $e$ is a negative loop at $v$, then $\Sigma/e$ has vertex set $V - \{ v \}$ and edge set resulting from $E$ by deleting $e$ and converting all edges incident with $v$ to half edges.

Before proving Theorem \ref{Spolynomial}, we give two illustrating examples.
First, a signed path with three vertices (Figure~\ref{fig:pathex}):

\begin{figure}[h]
\centering
\begin{tikzpicture}[scale=.5, line cap=round,line join=round,>=triangle 45,x=1.0cm,y=1.0cm]
\draw (0,0.5) node[above]{$\Sigma$};
\draw [dashed](-2.5,0)-- (-1,0);
\draw (-1,0)-- (0.5,0);
\draw (.5,0)-- (1.5,1);
\draw[color=black] (1.2,0.35) node {$e$};
\draw (-1.3,-0.3) node[below] {
$c_\Sigma(2k+1,2l)$};
\fill [color=black] (-2.5,0) circle (3pt);
\fill [color=black] (-1,0) circle (3pt);
\fill [color=black] (.5,0) circle (3pt);

\draw (2.5,1)--(4,1);
\draw (2.5,0.75)--(4,.75);

\draw [dashed](5,0)-- (6.5,0);
\draw (6.5,0)-- (8,0);
\fill [color=black] (5,0) circle (3pt);
\fill [color=black] (6.5,0) circle (3pt);
\fill [color=black] (8,0) circle (3pt);
\draw (7,-0.3) node[below] {
$c_{\Sigma-e}(2k+1,2l)$};

\draw (10,1)-- (11.5,1);

\draw [dashed](13.5,0)-- (15,0);
\draw (15,0)-- (16.5,0);
\fill [color=black] (13.5,0) circle (3pt);
\fill [color=black] (15,0) circle (3pt);
\draw (15,-0.3) node[below] {
$c_{\Sigma/e}(2k+1,2l)$};
\end{tikzpicture}
\caption{Deletion--contraction at a halfedge $e$ of a signed path $\Sigma$.}\label{fig:pathex}
\end{figure}
\vskip-.3cm
\begin{eqnarray}
c_\Sigma(2k+1,2l)&=&(2l)\big[(2l)(2l+2k+1)+(2k+1)(2l+2k)\big]+(2k)\big[(2l)(2l+2k+1)+(2k)(2l+2k)\big]\nonumber \\
&=&(2l)\big[(2l)(2l+2k+1)+(2k+1)(2l+2k)\big]+(2k+1)\big[(2l)(2l+2k+1)+(2k)(2l+2k)\big]\nonumber \\
&&- \;\big[(2l)(2l+2k+1)+(2k)(2l+2k)\big]\nonumber \\
&=& c_{\Sigma-e}(2k+1,2l) - c_{\Sigma/e}(2k+1,2l)\nonumber
\end{eqnarray}




\vskip.5cm

Our second example is a signed 3-cycle (Figure~\ref{fig:cycleex}):

\begin{figure}[h]
\centering
\begin{tikzpicture}[scale=.5, line cap=round,line join=round,>=triangle 45,x=1.0cm,y=1.0cm]
\draw (0,1) node[above]{$\Sigma$};
\draw [dashed](0,0)-- (1,2);
\draw [dashed](1,2)-- (2,0);
\draw (0,0)-- (2,0);
\draw[color=black] (1,0) node [above] {$e$};
\fill [color=black] (0,0) circle (3pt);
\fill [color=black] (1,2) circle (3pt);
\fill [color=black] (2,0) circle (3pt);
\draw (-1.5,-0.3) node[below] {
$c_\Sigma(2k+1,2l)$};

\draw (3,1)--(4.5,1);
\draw (3,0.75)--(4.5,0.75);

\draw [dashed](6,0)-- (7,2);
\draw [dashed](7,2)-- (8,0);
\fill [color=black] (6,0) circle (3pt);
\fill [color=black] (7,2) circle (3pt);
\fill [color=black] (8,0) circle (3pt);
\draw (6,-0.3) node[below] {
$c_{\Sigma-e}(2k+1,2l)$};

\draw (9,1)--(10.5,1);

\draw [dashed](13,0) to [out=120, in=230] (13,2);
\draw [dashed](13,0) to [out=60, in=-50] (13,2);
\draw (12.5,-0.3) node[below] {
$c_{\Sigma/e}(2k+1,2l)$};

\fill [color=black] (13,0) circle (3pt);
\fill [color=black] (13,2) circle (3pt);

\draw (15,1)--(16.5,1);
\draw (15.75,1.75)--(15.75,.25);

\draw (18,.5) node[above]{$2l$};
\fill [color=black] (19,1) circle (3pt);
\draw (20,-0.3) node[below] {
$2l\{c_{[\Sigma/e]-v}(2k+1,2l)\}$};

\end{tikzpicture}
\caption{Deletion--contraction at an edge $e$ of a signed 3-cycle $\Sigma$.}\label{fig:cycleex}
\end{figure}
\vskip-.3cm
\begin{eqnarray}
c_\Sigma(2k+1,2l)&=&(2l)\big[(2l)(2l+2k+1)+(2k+1)(2l+2k)\big]+(2k+1)\big[(2l)(2l+2k)+(2k)(2l+2k-1)\big]\nonumber \\
&=&(2l)\big[(2l)(2l+2k+1)+(2k+1)(2l+2k)\big]+(2k+1)\big[(2l)(2l+2k+1)+(2k)(2l+2k)\big]\nonumber \\
&&- \;\big[(2l)(2l+2k+1)+(2k+1)(2l+2k)\big]\nonumber \\
&& + \; 2l\big[2l+2k+1\big]\nonumber \\
&=& c_{\Sigma-e}(2k+1,2l) - c_{\Sigma/e}(2k+1,2l)+2l\{c_{[\Sigma/e]-v}(2k+1,2l)\} \nonumber
\end{eqnarray}

\vskip.3cm

\begin{proof}[Proof of Theorem \ref{Spolynomial}]
Let $e$ be a halfedge or negative loop of $\Sigma$, and let $v$ be its incident vertex.
Then $v$ cannot be colored zero, and so we have to subtract from the colorings of $\Sigma - e$ those which color $v$ zero (which are in bijection with the colorings of $\Sigma/e$).

Now let $e$ be an edge of $\Sigma$ that is not a halfedge or negative loop.
We have to subtract from the colorings of $\Sigma - e$ those which color the endpoints of $e$ the same (which are in bijection with the colorings of $\Sigma/e$) and add back in those where the latter color is $>k$ (the number of which is $2l$ times the number of proper $(k,l)$-colorings of $[\Sigma/e] - v$).
\end{proof}

By induction on the number of edges of a signed graph, we immediately conclude:

\begin{corollary}\label{chromaticpolcor}
The chromatic counting functions $c_\Sigma(2k+1,2l)$ and $c^*_\Sigma(2k,2l)$ are polynomials in $k$ and~$l$.
\end{corollary}

Our next result is the signed-graph analog of \cite[Theorem 1]{dohmenponitztittmann}.

\begin{theorem}\label{2variable SDPT}
Let $\Sigma=(V, E, \sigma)$ be a signed graph.  Then
\[
  c_\Sigma(2k+1,2l) = \sum_{W\subseteq V} (2l)^{|W|} \, c_{\Sigma - W}(2k+1,0)
  \qquad \text{ and } \qquad
  c^*_\Sigma(2k,2l) = \sum_{W\subseteq V} (2l)^{|W|} \, c^*_{\Sigma - W}(2k,0) \, .
\]
\end{theorem}

Naturally, the polynomials in the summations should be interpreted as Zaslavsky's chromatic polynomials.

\begin{proof}
Every proper $(k,l)$-coloring of $\Sigma$ can be obtained by first choosing a subset $W$ of $V$ that is colored with colors $>k$; there are $(2l)^{|W|}$ such colorings for these vertices. The remaining subgraph $\Sigma - W$ has to be colored properly using colors $\le k$.
\end{proof}

The above proof is virtually identical to that of \cite[Theorem 1]{dohmenponitztittmann}, and thus we obtain, as an analogous consequence, the following corollary, paralleling \cite[Corollary 2]{dohmenponitztittmann}, regarding the \emph{independence polynomial}
\[
  i_\Sigma (x) := \sum_{ W \subseteq V \text{ independent } } x^{ |V-W| } \, .
\]
(Here $W$ is \emph{independent} if no two vertices in $W$ are adjacent.)

\begin{corollary}\label{Sindependentpoly}
$i_\Sigma (x) = c_\Sigma(1,x) \, .$
\end{corollary}

\begin{proof}
By Theorem \ref{2variable SDPT}, we have
\[
  c_\Sigma(1,2x) = \sum_{W\subseteq V} (2x)^{|W|} \, c_{\Sigma - W}(1,0) \, .
\]
Now note that $c_{\Sigma - W}(1,0)$ equals one if $V - W$ is independent and zero otherwise.
\end{proof}

\begin{proof}[Proof of Theorem \ref{chromaticantibalancethm}]
By Theorem \ref{2variable SDPT},
\[
  c^*_\Sigma(2,2l)
  = \sum_{W\subseteq V} (2l)^{|W|} \, c^*_{\Sigma - W}(2,0) \, .
\]
Now $c^*_{S}(2,0)$ equals $2^{ c(S) }$ if $S$ is antibalanced, and $c^*_{S}(2,0) = 0$ if $S$ is not antibalanced. Thus
\[
  c^*_\Sigma(2,2l)
  = \sum_{ S \text{ antibalanced subgraph of } \Sigma } (2l)^{ |V| - |V(S)| } 2^{ c(S) }
  = (2l)^{ |V| } \, a_\Sigma ( \tfrac{ 1 }{ 2l } , 2) \, . \qedhere
\]
\end{proof}

For completeness sake, we state the signed analogue of \cite[Theorem 3]{dohmenponitztittmann}; its proof is virtually
identical to the unsigned case.

\begin{theorem}
$\displaystyle \frac{ \partial }{ \partial (2l) } \, c_\Sigma(2k+1, 2l) = \sum_{ v \in V } c_{ \Sigma - v } (2k+1, 2l)$
\ and \
$\displaystyle \frac{ \partial }{ \partial (2l) } \, c^*_\Sigma(2k, 2l) = \sum_{ v \in V } c^*_{ \Sigma - v } (2k, 2l) \, .$
\end{theorem}


\section{Bivariate Chromatic Reciprocity Theorems}\label{reciprocitysection}

The proofs of our reciprocity theorems follow along the lines of the proof of Stanley's Theorem \ref{stanleyrecthm} given in \cite{iop}, which introduced the general setup of an \emph{inside-out polytope} $(\P, \HH)$ consisting of a
rational polytope $\P$ and a rational hyperplane arrangement $\HH$ in
$\R^d$; that is, the linear equations and inequalities defining $\P$ and $\HH$ have integer coefficients.
(The proper understanding of this section assumes familiarity with \cite{iop}.)
The goal is to compute the counting function
\[
  L_{ \P, \HH }^\circ (t) := \left| (\P \setminus \HH) \cap \tfrac 1 t \Z^d \right| ,
\]
and it follows from Ehrhart's theory of counting lattice points in
dilates of rational polytopes \cite{ccd,ehrhartpolynomial} that this
function is a quasipolynomial in $t$ whose degree is $\dim(P)$, whose
(constant) leading coefficient is the normalized lattice volume of $P$,
and whose period divides the lcm of all denominators that appear in
the coordinates of the vertices of $(\P, \HH)$. 
In our case, all vertices of $(\P, \HH)$ will be integral, so that the resulting counting functions will be polynomials.
Furthermore, \cite{iop} established the reciprocity theorem
\begin{equation}\label{iopreciprocity}
  L_{ \P^\circ, \HH }^\circ (-t) = (-1)^{ \dim \P } L_{ \overline \P, \HH } (t) \, ,
\end{equation}
where $\P^\circ$ and $\overline \P$ denote the interior and closure of $\P$, respectively, and
\begin{equation}\label{multdef}
  L_{ \P, \HH } (t) := \sum_{ \m \in \frac 1 t \Z^d } \mult_{ \P, \HH } (\m)
\end{equation}
where $\mult_{ \P, \HH } (\m)$ denotes the number of closed regions of $(\P, \HH)$ containing $\m$.
(A \emph{region} of $(\P, \HH)$ is a connected component of $\P \setminus \HH$; a \emph{closed region} is the closure of a region.)
See \cite{iop} for this and several more properties of inside-out
polytopes. The concept of inside-out polytopes has been applied to a number of combinatorial settings; at the heart of any such application is an interpretation of the regions of $(\P, \HH)$; from this point of view, Stanley's Theorem \ref{stanleyrecthm} follows from Greene's observation \cite{greeneacyclic,greenezaslavsky} that the regions of the \emph{graphic arrangement} (for a given graph $\Gamma = (V, E)$)
\[
  \HH_\Gamma := \left\{ x_v = x_w : \, vw \in E \right\}
\]
in $\R^V$ are in one-to-one correspondence with the acyclic orientations of $\Gamma$.

\begin{proof}[Proof of Theorem \ref{bivariaterecthm}]
Given $\Gamma = (V, E)$, let $\Box := [0,1]^V$ be the unit cube in $\R^V$ and
\[
  \Phi(k,l) := \left| (k+l+1) \Box^\circ \cap \Z^V \right| - \left| (k+1) \Box^\circ \cap \Z^V \right| ,
\]
i.e., $\Phi(k,l)$ is the difference of two evaluations (at $k+l+1$ and $k+1$) of the Ehrhart polynomial of $\Box^\circ$.
By Ehrhart--Macdonald reciprocity (see, e.g., \cite[Chapter 4]{ccd}),
\begin{equation}\label{phireceq}
  (-1)^{ |V| } \Phi(-k,-l) = \left| (k+l-1) \Box \cap \Z^V \right| - \left| (k-1) \Box \cap \Z^V \right| .
\end{equation}
On the other hand, it is natural to interpret the bichromatic polynomial of $\Gamma$ geometrically as
\begin{align*}
  c_\Gamma (k, l)
  &= \left| \left( [1,k+l]^V - \left( [1,k]^V \cap \HH_\Gamma \right) \right) \cap \Z^V \right|
  = \left| \left( (0,k+l+1)^V - \left( (0,k+1)^V \cap\ \HH_\Gamma \right) \right) \cap \Z^V \right| \\
  &= L_{ \Box^\circ, \HH_\Gamma }^\circ (k+1) + \Phi(k,l)
\end{align*}
(see Figure \ref{polytopalG} for an illustrative example).
\begin{figure}[h]
\centering
\begin{tikzpicture}[scale=.5, line cap=round,line join=round,>=triangle 45,x=1.0cm,y=1.0cm]
\draw[->,color=black,dashed] (-3.54,0) -- (12.2,0);
\draw[->,color=black,dashed] (0,-2.9) -- (0,12.2);
\draw [->] (7,-2.96) -- (7,-0.96);
\draw [->] (-2.98,7) -- (-0.98,7);
\draw (-5,7.6) node[anchor=north west] {\footnotesize$k+1$};
\draw (6.6,-2.84) node[anchor=north west] {\footnotesize$k+1$};
\draw [->] (11,-3) -- (11,-1);
\draw [->] (-3,11) -- (-1,11);
\draw (10.52,-2.84) node[anchor=north west] {\footnotesize$k+l+1$};
\draw (-6,11.6) node[anchor=north west] {\footnotesize$k+l+1$};

\fill [color=black] (1,2) circle (4pt);\fill [color=black] (1,3) circle (4pt);
\fill [color=black] (1,4) circle (4pt);\fill[color=black] (1,5) circle (4pt);
\fill [color=black] (1,6) circle (4pt);\fill [color=black] (2,6) circle (4pt);
\fill [color=black] (3,6) circle (4pt);\fill [color=black] (2,5) circle (4pt);
\fill [color=black] (2,4) circle (4pt);\fill [color=black] (2,3) circle (4pt);
\fill [color=black] (3,2) circle (4pt);\fill [color=black] (3,4) circle (4pt);
\fill [color=black] (3,5) circle (4pt);\fill [color=black] (4,6) circle (4pt);
\fill [color=black] (5,6) circle (4pt);\fill [color=black] (4,5) circle (4pt);
\fill [color=black] (2,1) circle (4pt);\fill [color=black] (3,1) circle (4pt);
\fill [color=black] (4,1) circle (4pt);\fill [color=black] (4,2) circle (4pt);
\fill [color=black] (4,3) circle (4pt);\fill [color=black] (5,4) circle (4pt);
\fill [color=black] (5,3) circle (4pt);\fill [color=black] (5,2) circle (4pt);
\fill [color=black] (5,1) circle (4pt);\fill [color=black] (6,1) circle (4pt);
\fill [color=black] (6,2) circle (4pt);\fill [color=black] (6,3) circle (4pt);
\fill [color=black] (6,4) circle (4pt);\fill [color=black] (6,5) circle (4pt);

\fill [color=black] (1,7) circle (4pt);\fill [color=black] (2,7) circle (4pt);
\fill [color=black] (3,7) circle (4pt);\fill [color=black] (4,7) circle (4pt);
\fill [color=black] (5,7) circle (4pt);\fill [color=black] (6,7) circle (4pt);
\fill [color=black] (7,7) circle (4pt);\fill [color=black] (7,6) circle (4pt);
\fill [color=black] (7,5) circle (4pt);\fill [color=black] (7,4) circle (4pt);
\fill [color=black] (7,3) circle (4pt);\fill [color=black] (7,2) circle (4pt);
\fill [color=black] (7,1) circle (4pt);\fill [color=black] (1,8) circle (4pt);
\fill [color=black] (1,9) circle (4pt);\fill [color=black] (1,10) circle (4pt);
\fill [color=black] (2,10) circle (4pt);\fill [color=black] (2,9) circle (4pt);
\fill [color=black] (2,8) circle (4pt);\fill [color=black] (3,8) circle (4pt);
\fill [color=black] (3,9) circle (4pt);\fill [color=black] (3,10) circle (4pt);
\fill [color=black] (4,10) circle (4pt);\fill [color=black] (4,9) circle (4pt);
\fill [color=black] (4,8) circle (4pt);\fill [color=black] (5,8) circle (4pt);
\fill [color=black] (5,9) circle (4pt);\fill [color=black] (5,10) circle (4pt);
\fill [color=black] (6,10) circle (4pt);\fill [color=black] (6,9) circle (4pt);
\fill [color=black] (6,8) circle (4pt);\fill [color=black] (7,8) circle (4pt);
\fill [color=black] (7,9) circle (4pt);\fill [color=black] (7,10) circle (4pt);
\fill [color=black] (8,10) circle (4pt);\fill [color=black] (9,10) circle (4pt);
\fill [color=black] (10,10) circle (4pt);\fill [color=black] (10,9) circle (4pt);
\fill [color=black] (9,9) circle (4pt);\fill [color=black] (8,9) circle (4pt);
\fill [color=black] (8,8) circle (4pt);\fill [color=black] (8,7) circle (4pt);
\fill [color=black] (9,7) circle (4pt);\fill [color=black] (9,8) circle (4pt);
\fill [color=black] (10,8) circle (4pt);\fill [color=black] (10,7) circle (4pt);
\fill [color=black] (10,6) circle (4pt);\fill [color=black] (10,5) circle (4pt);
\fill [color=black] (10,4) circle (4pt);\fill [color=black] (10,3) circle (4pt);
\fill [color=black] (10,2) circle (4pt);\fill [color=black] (10,1) circle (4pt);
\fill [color=black] (9,2) circle (4pt);\fill [color=black] (9,3) circle (4pt);
\fill [color=black] (9,5) circle (4pt);\fill [color=black] (9,4) circle (4pt);
\fill [color=black] (9,1) circle (4pt);\fill [color=black] (9,6) circle (4pt);
\fill [color=black] (8,6) circle (4pt);\fill [color=black] (8,5) circle (4pt);
\fill [color=black] (8,4) circle (4pt);\fill [color=black] (8,3) circle (4pt);
\fill [color=black] (8,2) circle (4pt);\fill [color=black] (8,1) circle (4pt);

\draw [color=black] (0,11) circle (4pt);\draw [color=black] (0,10) circle (4pt);
\draw [color=black] (0,9) circle (4pt);\draw [color=black] (0,8) circle (4pt);
\draw [color=black] (0,7) circle (4pt);\draw [color=black] (0,6) circle (4pt);
\draw [color=black] (0,5) circle (4pt);\draw [color=black] (0,4) circle (4pt);
\draw [color=black] (0,3) circle (4pt);\draw [color=black] (0,2) circle (4pt);
\draw [color=black] (0,1) circle (4pt);\draw [color=black] (0,0) circle (4pt);
\draw [color=black] (1,0) circle (4pt);\draw [color=black] (2,0) circle (4pt);
\draw [color=black] (3,0) circle (4pt);\draw [color=black] (4,0) circle (4pt);
\draw [color=black] (5,0) circle (4pt);\draw [color=black] (6,0) circle (4pt);
\draw [color=black] (7,0) circle (4pt);\draw [color=black] (8,0) circle (4pt);
\draw [color=black] (9,0) circle (4pt);\draw [color=black] (10,0) circle (4pt);
\draw [color=black] (11,0) circle (4pt);

\draw [color=black] (11,1) circle (4pt);\draw [color=black] (11,2) circle (4pt);
\draw [color=black] (11,3) circle (4pt);\draw [color=black] (11,4) circle (4pt);
\draw [color=black] (11,5) circle (4pt);\draw [color=black] (11,6) circle (4pt);
\draw [color=black] (11,7) circle (4pt);\draw [color=black] (11,8) circle (4pt);
\draw [color=black] (11,9) circle (4pt);\draw [color=black] (11,10) circle (4pt);
\draw [color=black] (11,11) circle (4pt);\draw [color=black] (10,11) circle (4pt);
\draw [color=black] (9,11) circle (4pt);\draw [color=black] (8,11) circle (4pt);
\draw [color=black] (7,11) circle (4pt);\draw [color=black] (6,11) circle (4pt);
\draw [color=black] (5,11) circle (4pt);\draw [color=black] (4,11) circle (4pt);
\draw [color=black] (3,11) circle (4pt);\draw [color=black] (2,11) circle (4pt);
\draw [color=black] (1,11) circle (4pt);\draw [color=black] (1,1) circle (4pt);
\draw [color=black] (2,2) circle (4pt);\draw [color=black] (3,3) circle (4pt);
\draw [color=black] (4,4) circle (4pt);\draw [color=black] (5,5) circle (4pt);
\draw [color=black] (6,6) circle (4pt);
\draw (0.1,7)--(7,7);\draw (7,7)--(7,0.1);
\draw [dashed](0.1,11)--(11,11);\draw [dashed](11,11)--(11,0.1);
\end{tikzpicture}
\caption{The proper $(k,l)$-colorings of $K_2$ with $k=6$ and $l=4$.}
\label{polytopalG}
\end{figure}
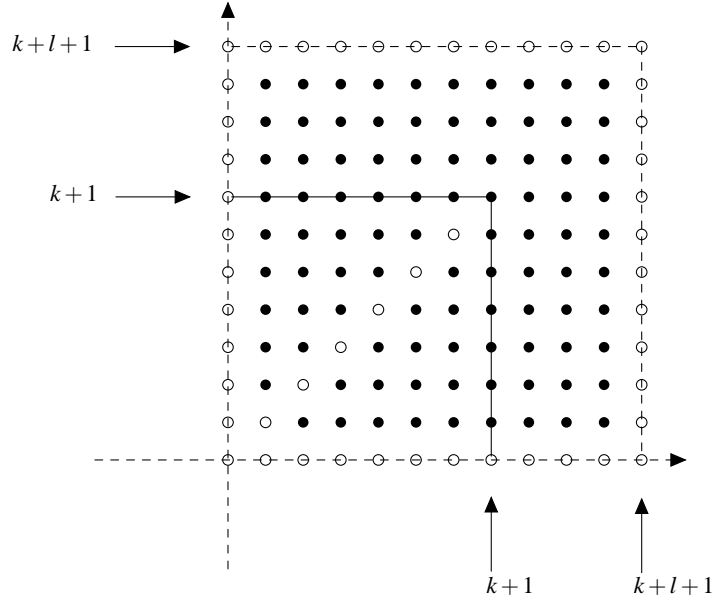
Thus, by \eqref{iopreciprocity} and \eqref{phireceq},
\begin{align*}
  (-1)^{ |V| } c_\Gamma (-k, -l)
  &= (-1)^{ |V| } L_{ \Box^\circ, \HH_\Gamma }^\circ (-k+1) + (-1)^{ |V| } \Phi(-k,-l) \\
  &= L_{ \Box, \HH_\Gamma } (k-1) + \left| (k+l-1) \Box \cap \Z^V \right| - \left| (k-1) \Box \cap \Z^V \right| .
\end{align*}
What we are counting on the right-hand side are the integer lattice points in the cube $[0,k+l-1]^V$, with multiplicity equal to $1$ if outside the cube $[0,k-1]^V$, otherwise with multiplicity equal to the number of closed regions of $\HH_\Gamma$ the points lies in. As we mentioned above (and as was used in \cite{iop}), the latter can be interpreted as the number of compatible acyclic orientations of $\Gamma$. It is now a short step to re-interpret the lattice points in $[0,k+l-1]^V$ as $(k+l)$-colorings and the ones in $[0,k-1]$ as $k$-colorings of~$\Gamma$.
\end{proof}

In order to state and prove the analogous reciprocity theorem for bichromatic polynomials for signed graphs, we need more definitions.
An \emph{orientation} of a signed graph $\Sigma = (\Gamma, \sigma)$ is obtained from a \emph{bidirection} of the underlying graph $\Gamma$, where the endpoints of each edge are independently oriented, in such a way that the two arrows on an edge $e$ point in the same direction if $\sigma_e = +$ and they conflict if $\sigma_e = -$.
We express the bidirection (and hence the orientation) by means of an \emph{incidence function} $\eta$ defined on the edge ends: $\eta_{ ve } = 1$ if the arrow on $e$ at $v$ points into $v$, and $\eta_{ ve } = -1$ if it points away from $v$; with this definition we obtain $\sigma_e = - \eta_{ ve } \eta_{ we }$ for an edge $e = vw$. (See \cite{zaslavskyorientationsignedgraphs}  for more details.)

Following \cite{zaslavskysignedcoloring}, we call a coloring $\x \in \Z^V$ and an orientation $\eta$ \emph{compatible} if for any link $e = vw$
\[
  \eta_{ ve } x_v + \eta_{ we } x_w \ge 0 \, ,
\]
and for any halfedge or negative loop $e$ at $v$
\[
  \eta_{ ve } x_v \ge 0 \, .
\]
Furthermore, an orientation is \emph{acyclic} if no cycle has a source or sink (i.e., a vertex $v$ for which both incident edges point either into or away from $v$).
Zaslavsky \cite{zaslavskysignedcoloring} proved the following analogue of Stanley's Theorem~\ref{stanleyrecthm} for the chromatic polynomial $c_\Sigma(2k+1)$ of a signed graph~$\Sigma$:

\begin{theorem}[Zaslavsky]
For $k \in \Z_{ >0 }$,
$(-1)^{ |V| } c_\Sigma (-2k-1)$ equals the number of $k$-colorings of $\Sigma$, each counted with multiplicity equal to the number of compatible acyclic orientations of $\Sigma$.
In particular, $(-1)^{ |V| } c_\Sigma (-1)$ equals the number of acyclic orientations of~$\Sigma$.
\end{theorem}

Our analogue of this theorem for bivariate chromatic polynomials is as follows.

\begin{theorem}
For $k \in \Z_{ >0 }$ and $l \in \Z_{ \ge 0 }$,
$(-1)^{ |V| } c_\Sigma (-2k-1,-2l)$ equals the number of $(k+l)$-colorings of $\Sigma$, each counted with multiplicity: a $k$-coloring has multiplicity equal to the number of compatible acyclic orientations of $\Sigma$, and a coloring that uses at least one color with absolute value $> k$ has multiplicity~$1$.
\end{theorem}

\begin{proof}
We follow the proof of Theorem \ref{bivariaterecthm}, with a few modifications:
Given $\Sigma = (V, E, \sigma)$, let $\Box := [-1,1]^V$ and
\[
  \Phi(k,l) := \left| (k+l+1) \Box^\circ \cap \Z^V \right| - \left| (k+1) \Box^\circ \cap \Z^V \right| .
\]
By Ehrhart--Macdonald reciprocity,
\begin{equation}\label{phirecsignedeq}
  (-1)^{ |V| } \Phi(-k,-l) = \left| (k+l-1) \Box \cap \Z^V \right| - \left| (k-1) \Box \cap \Z^V \right| .
\end{equation}
To construct an inside-out counting function for $\Sigma$, we use the hyperplane arrangement
\[
  \HH_\Sigma := \left\{ x_v = \sigma_{ vw } \, x_w : \, vw \in E \right\} ,
\]
and so
\begin{align*}
  c_\Sigma (2k+1, 2l)
  &= \left| \left( [-k-l,k+l]^V - \left( [-k,k]^V \cap \HH_\Sigma \right) \right) \cap \Z^V \right| \\
  &= \left| \left( (-k-l-1,k+l+1)^V - \left( (-k-1,k+1)^V \cap\ \HH_\Sigma \right) \right) \cap \Z^V \right| \\
  &= L_{ \Box^\circ, \HH_\Sigma }^\circ (k+1) + \Phi(k,l)
\end{align*}
(see Figure \ref{polytopalS}).
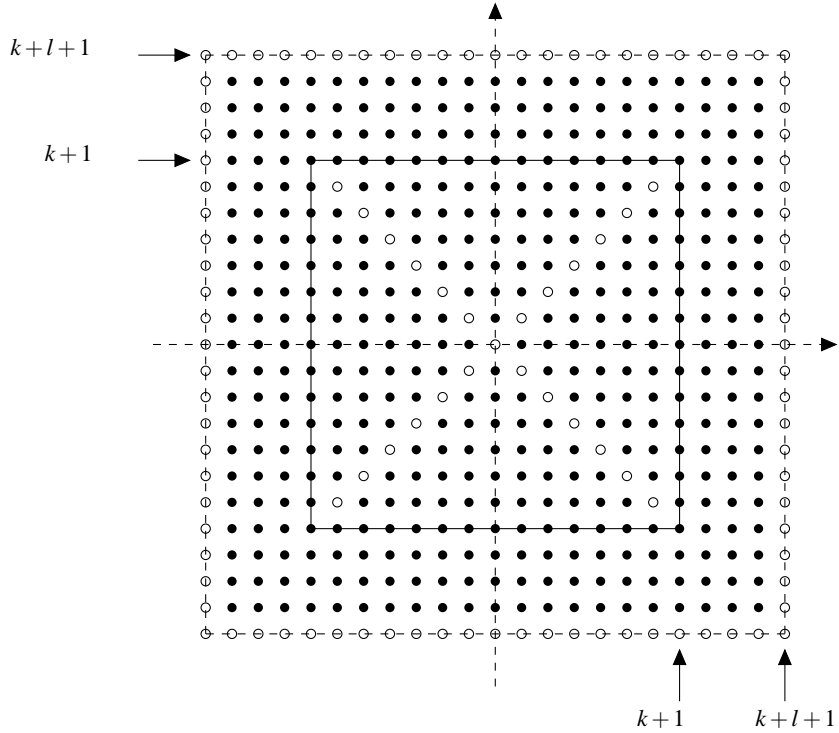
\begin{figure}[h]
\centering
\begin{tikzpicture}[scale=.35,>=triangle 45]
\draw[->,color=black,dashed] (-13,0) -- (13,0);
\draw[->,color=black,dashed] (0,-13) -- (0,13);
\draw [->] (11,-13.56)--(11,-11.56);
\draw [->] (7,-13.56)--(7,-11.56);
\draw [->] (-13.56,11)--(-11.56,11);
\draw [->] (-13.56,7)--(-11.56,7);
\draw (5,-13.5) node[anchor=north west] {\footnotesize$k+1$};
\draw (9.5,-13.5) node[anchor=north west] {\footnotesize$k+l+1$};
\draw (-17.5,8) node[anchor=north west] {\footnotesize$k+1$};
\draw (-18.8,12) node[anchor=north west] {\footnotesize$k+l+1$};

\draw (-7,7)-- (7,7);\draw (7,7)-- (7,-7);
\draw (7,-7)-- (-7,-7);\draw (-7,-7)-- (-7,7);
\draw [dashed](-11,11)-- (11,11);\draw [dashed](11,11)-- (11,-11);
\draw [dashed](11,-11)-- (-11,-11);\draw [dashed](-11,-11)-- (-11,11);
\fill [color=black] (1,7) circle (5pt);\fill [color=black] (2,7) circle (5pt);
\fill [color=black] (3,7) circle (5pt);\fill [color=black] (4,7) circle (5pt);
\fill [color=black] (5,7) circle (5pt);\fill [color=black] (6,7) circle (5pt);
\fill [color=black] (7,6) circle (5pt);\fill [color=black] (7,5) circle (5pt);
\fill [color=black] (7,4) circle (5pt);\fill [color=black] (7,3) circle (5pt);
\fill [color=black] (7,2) circle (5pt);\fill [color=black] (7,1) circle (5pt);
\fill [color=black] (1,8) circle (5pt);\fill [color=black] (1,9) circle (5pt);
\fill [color=black] (1,10) circle (5pt);\fill [color=black] (2,10) circle (5pt);
\fill [color=black] (2,9) circle (5pt);\fill [color=black] (2,8) circle (5pt);
\fill [color=black] (3,8) circle (5pt);\fill [color=black] (3,9) circle (5pt);
\fill [color=black] (3,10) circle (5pt);\fill [color=black] (4,10) circle (5pt);
\fill [color=black] (4,9) circle (5pt);\fill [color=black] (4,8) circle (5pt);
\fill [color=black] (5,8) circle (5pt);\fill [color=black] (5,9) circle (5pt);
\fill [color=black] (5,10) circle (5pt);\fill [color=black] (6,10) circle (5pt);
\fill [color=black] (6,9) circle (5pt);\fill [color=black] (6,8) circle (5pt);
\fill [color=black] (7,8) circle (5pt);\fill [color=black] (7,9) circle (5pt);
\fill [color=black] (7,10) circle (5pt);\fill [color=black] (8,10) circle (5pt);
\fill [color=black] (9,10) circle (5pt);\fill [color=black] (10,9) circle (5pt);
\fill [color=black] (8,9) circle (5pt);\fill [color=black] (8,7) circle (5pt);
\fill [color=black] (9,7) circle (5pt);\fill [color=black] (9,8) circle (5pt);
\fill [color=black] (10,8) circle (5pt);\fill [color=black] (10,7) circle (5pt);
\fill [color=black] (10,6) circle (5pt);\fill [color=black] (10,5) circle (5pt);
\fill [color=black] (10,4) circle (5pt);\fill [color=black] (10,3) circle (5pt);
\fill [color=black] (10,2) circle (5pt);\fill [color=black] (10,1) circle (5pt);
\fill [color=black] (9,2) circle (5pt);\fill [color=black] (9,3) circle (5pt);
\fill [color=black] (9,5) circle (5pt);\fill [color=black] (9,4) circle (5pt);
\fill [color=black] (9,1) circle (5pt);\fill [color=black] (9,6) circle (5pt);
\fill [color=black] (8,6) circle (5pt);\fill [color=black] (8,5) circle (5pt);
\fill [color=black] (8,4) circle (5pt);\fill [color=black] (8,3) circle (5pt);
\fill [color=black] (8,2) circle (5pt);\fill [color=black] (8,1) circle (5pt);
\draw (0,11) circle (5pt);
\fill (0,10) circle (5pt);\fill (0,9) circle (5pt);
\fill (0,8) circle (5pt);\fill (0,7) circle (5pt);

\fill (0,6) circle (5pt);\fill (0,5) circle (5pt);
\fill (0,4) circle (5pt);\fill (0,3) circle (5pt);
\fill (0,2) circle (5pt);\fill(0,1) circle (5pt);
\fill (1,0) circle (5pt);\fill (2,0) circle (5pt);
\fill (3,0) circle (5pt);\fill (4,0) circle (5pt);
\fill (5,0) circle (5pt);\fill (6,0) circle (5pt);
\fill (7,0) circle (5pt);\fill (8,0) circle (5pt);
\fill (9,0) circle (5pt);\fill (10,0) circle (5pt);
\draw (11,0) circle (5pt);\draw (11,1) circle (5pt);
\draw (11,2) circle (5pt);\draw (11,3) circle (5pt);
\draw (11,4) circle (5pt);\draw (11,5) circle (5pt);
\draw (11,6) circle (5pt);\draw (11,7) circle (5pt);
\draw (11,8) circle (5pt);\draw (11,9) circle (5pt);
\draw (11,10) circle (5pt);\draw (10,11) circle (5pt);
\draw (9,11) circle (5pt);\draw (8,11) circle (5pt);
\draw (7,11) circle (5pt);\draw (6,11) circle (5pt);
\draw (5,11) circle (5pt);\draw (4,11) circle (5pt);
\draw (3,11) circle (5pt);\draw (2,11) circle (5pt);
\draw (1,11) circle (5pt);
\fill [color=black] (7,7) circle (5pt);\fill [color=black] (8,8) circle (5pt);
\fill [color=black] (9,9) circle (5pt);\fill [color=black] (10,10) circle (5pt);
\draw (11,11) circle (5pt);\fill (1,6) circle (5pt);
\fill (2,6) circle (5pt);\fill (3,6) circle (5pt);
\fill(4,6) circle (5pt);\fill (5,6) circle (5pt);
\draw (6,6) circle (5pt);\fill (6,5) circle (5pt);
\draw (5,5) circle (5pt);\fill (4,5) circle (5pt);
\fill (3,5) circle (5pt);\fill (2,5) circle (5pt);
\fill (1,5) circle (5pt);\fill (1,4) circle (5pt);
\fill (1,3) circle (5pt);\fill (1,2) circle (5pt);
\draw (1,1) circle (5pt);\draw (0,0) circle (5pt);
\fill (2,1) circle (5pt);\fill (3,1) circle (5pt);
\fill (4,1) circle (5pt);\fill (5,1) circle (5pt);
\fill (6,1) circle (5pt);\fill (6,2) circle (5pt);
\fill (6,3) circle (5pt);\fill (6,4) circle (5pt);
\fill (5,4) circle (5pt);\draw (4,4) circle (5pt);
\fill (3,4) circle (5pt);\fill (2,4) circle (5pt);
\fill (2,3) circle (5pt);\draw (3,3) circle (5pt);
\fill (4,3) circle (5pt);\fill (5,3) circle (5pt);
\fill (5,2) circle (5pt);\fill (4,2) circle (5pt);
\fill (3,2) circle (5pt);\draw (2,2) circle (5pt);
\fill [color=black] (-1,7) circle (5pt);\fill [color=black] (-2,7) circle (5pt);
\fill [color=black] (-3,7) circle (5pt);\fill [color=black] (-4,7) circle (5pt);
\fill [color=black] (-5,7) circle (5pt);\fill [color=black] (-6,7) circle (5pt);
\fill [color=black] (-7,7) circle (5pt);\fill [color=black] (-8,7) circle (5pt);
\fill [color=black] (-9,7) circle (5pt);\fill [color=black] (-9,6) circle (5pt);
\fill [color=black] (-9,5) circle (5pt);\fill [color=black] (-9,4) circle (5pt);
\fill [color=black] (-9,3) circle (5pt);\fill [color=black] (-9,2) circle (5pt);
\fill [color=black] (-9,1) circle (5pt);\fill [color=black] (-1,8) circle (5pt);
\fill [color=black] (-2,8) circle (5pt);\fill [color=black] (-3,8) circle (5pt);
\fill [color=black] (-4,8) circle (5pt);\fill [color=black] (-5,8) circle (5pt);
\fill [color=black] (-6,8) circle (5pt);\fill [color=black] (-7,8) circle (5pt);
\fill [color=black] (-8,8) circle (5pt);\fill [color=black] (-9,8) circle (5pt);
\fill [color=black] (-1,9) circle (5pt);\fill [color=black] (-2,9) circle (5pt);
\fill [color=black] (-3,9) circle (5pt);\fill [color=black] (-4,9) circle (5pt);
\fill [color=black] (-5,9) circle (5pt);\fill [color=black] (-6,9) circle (5pt);
\fill [color=black] (-7,9) circle (5pt);\fill [color=black] (-8,9) circle (5pt);
\fill [color=black] (-9,9) circle (5pt);\fill [color=black] (-1,10) circle (5pt);
\fill [color=black] (-2,10) circle (5pt);\fill [color=black] (-3,10) circle (5pt);
\fill [color=black] (-4,10) circle (5pt);\fill [color=black] (-5,10) circle (5pt);
\fill [color=black] (-6,10) circle (5pt);\fill [color=black] (-7,10) circle (5pt);
\fill [color=black] (-8,10) circle (5pt);\fill [color=black] (-9,10) circle (5pt);
\fill [color=black] (-10,10) circle (5pt);\fill [color=black] (-10,9) circle (5pt);
\fill [color=black] (-10,8) circle (5pt);\fill [color=black] (-10,7) circle (5pt);
\fill [color=black] (-10,6) circle (5pt);\fill [color=black] (-10,5) circle (5pt);
\fill [color=black] (-10,4) circle (5pt);\fill [color=black] (-10,3) circle (5pt);
\fill [color=black] (-10,2) circle (5pt);\fill [color=black] (-10,1) circle (5pt);
\fill [color=black] (-8,1) circle (5pt);\fill [color=black] (-8,2) circle (5pt);
\fill [color=black] (-8,3) circle (5pt);\fill [color=black] (-8,4) circle (5pt);
\fill [color=black] (-8,5) circle (5pt);\fill [color=black] (-8,6) circle (5pt);
\fill [color=black] (-7,6) circle (5pt);\fill [color=black] (-7,5) circle (5pt);
\fill [color=black] (-7,4) circle (5pt);\fill [color=black] (-7,3) circle (5pt);
\fill [color=black] (-7,1) circle (5pt);\fill [color=black] (-10,-1) circle (5pt);
\fill [color=black] (-9,-1) circle (5pt);\fill [color=black] (-8,-1) circle (5pt);
\fill [color=black] (-7,-1) circle (5pt);\fill [color=black] (-10,-2) circle (5pt);
\fill [color=black] (-9,-2) circle (5pt);\fill [color=black] (-8,-2) circle (5pt);
\fill [color=black] (-7,-2) circle (5pt);\fill [color=black] (-10,-3) circle (5pt);
\fill [color=black] (-9,-3) circle (5pt);\fill [color=black] (-8,-3) circle (5pt);
\fill [color=black] (-7,-3) circle (5pt);\fill [color=black] (-10,-4) circle (5pt);
\fill [color=black] (-9,-4) circle (5pt);\fill [color=black] (-8,-4) circle (5pt);
\fill [color=black] (-7,-4) circle (5pt);\fill [color=black] (-10,-5) circle (5pt);
\fill [color=black] (-9,-5) circle (5pt);\fill [color=black] (-8,-5) circle (5pt);
\fill [color=black] (-7,-5) circle (5pt);\fill [color=black] (-10,-6) circle (5pt);
\fill [color=black] (-9,-6) circle (5pt);\fill [color=black] (-8,-6) circle (5pt);
\fill [color=black] (-7,-6) circle (5pt);\fill [color=black] (-10,-7) circle (5pt);
\fill [color=black] (-9,-7) circle (5pt);\fill [color=black] (-8,-7) circle (5pt);
\fill [color=black] (-7,-7) circle (5pt);\fill [color=black] (-10,-8) circle (5pt);
\fill [color=black] (-9,-8) circle (5pt);\fill [color=black] (-8,-8) circle (5pt);
\fill [color=black] (-7,-8) circle (5pt);\fill [color=black] (-10,-9) circle (5pt);
\fill [color=black] (-9,-9) circle (5pt);\fill [color=black] (-8,-9) circle (5pt);
\fill [color=black] (-7,-9) circle (5pt);\fill [color=black] (-10,-10) circle (5pt);
\fill [color=black] (-9,-10) circle (5pt);\fill [color=black] (-8,-10) circle (5pt);
\fill [color=black] (-7,-10) circle (5pt);\fill [color=black] (-6,-10) circle (5pt);
\fill [color=black] (-5,-10) circle (5pt);\fill [color=black] (-4,-10) circle (5pt);
\fill [color=black] (-3,-10) circle (5pt);\fill [color=black] (-2,-10) circle (5pt);
\fill [color=black] (-1,-10) circle (5pt);\fill [color=black] (-6,-9) circle (5pt);
\fill [color=black] (-5,-9) circle (5pt);\fill [color=black] (-4,-9) circle (5pt);
\fill [color=black] (-3,-9) circle (5pt);\fill [color=black] (-2,-9) circle (5pt);
\fill [color=black] (-1,-9) circle (5pt);\fill [color=black] (-6,-8) circle (5pt);
\fill [color=black] (-5,-8) circle (5pt);\fill [color=black] (-4,-8) circle (5pt);
\fill [color=black] (-3,-8) circle (5pt);\fill [color=black] (-2,-8) circle (5pt);
\fill [color=black] (-1,-8) circle (5pt);\fill [color=black] (-6,-7) circle (5pt);
\fill [color=black] (-5,-7) circle (5pt);\fill [color=black] (-4,-7) circle (5pt);
\fill [color=black] (-3,-7) circle (5pt);\fill [color=black] (-2,-7) circle (5pt);
\fill [color=black] (-1,-7) circle (5pt);\fill [color=black] (1,-7) circle (5pt);
\fill [color=black] (2,-7) circle (5pt);\fill [color=black] (3,-7) circle (5pt);
\fill [color=black] (4,-7) circle (5pt);\fill [color=black] (5,-7) circle (5pt);
\fill [color=black] (6,-7) circle (5pt);\fill [color=black] (7,-7) circle (5pt);
\fill [color=black] (8,-7) circle (5pt);\fill [color=black] (9,-7) circle (5pt);
\fill [color=black] (10,-7) circle (5pt);\fill [color=black] (1,-8) circle (5pt);
\fill [color=black] (2,-8) circle (5pt);\fill [color=black] (3,-8) circle (5pt);
\fill [color=black] (4,-8) circle (5pt);\fill [color=black] (5,-8) circle (5pt);
\fill [color=black] (6,-8) circle (5pt);\fill [color=black] (7,-8) circle (5pt);
\fill [color=black] (8,-8) circle (5pt);\fill [color=black] (9,-8) circle (5pt);
\fill [color=black] (10,-8) circle (5pt);\fill [color=black] (1,-9) circle (5pt);
\fill [color=black] (2,-9) circle (5pt);\fill [color=black] (3,-9) circle (5pt);
\fill [color=black] (4,-9) circle (5pt);\fill [color=black] (5,-9) circle (5pt);
\fill [color=black] (6,-9) circle (5pt);\fill [color=black] (7,-9) circle (5pt);
\fill [color=black] (8,-9) circle (5pt);\fill [color=black] (9,-9) circle (5pt);
\fill [color=black] (10,-9) circle (5pt);\fill [color=black] (1,-10) circle (5pt);
\fill [color=black] (2,-10) circle (5pt);\fill [color=black] (3,-10) circle (5pt);
\fill [color=black] (4,-10) circle (5pt);\fill [color=black] (5,-10) circle (5pt);
\fill [color=black] (6,-10) circle (5pt);\fill [color=black] (7,-10) circle (5pt);
\fill [color=black] (8,-10) circle (5pt);\fill [color=black] (9,-10) circle (5pt);
\fill [color=black] (10,-10) circle (5pt);\fill [color=black] (7,-1) circle (5pt);
\fill [color=black] (8,-1) circle (5pt);\fill [color=black] (9,-1) circle (5pt);
\fill [color=black] (10,-1) circle (5pt);\fill [color=black] (7,-2) circle (5pt);
\fill [color=black] (8,-2) circle (5pt);\fill [color=black] (9,-2) circle (5pt);
\fill [color=black] (10,-2) circle (5pt);\fill [color=black] (7,-3) circle (5pt);
\fill [color=black] (8,-3) circle (5pt);\fill [color=black] (9,-3) circle (5pt);
\fill [color=black] (10,-3) circle (5pt);\fill [color=black] (7,-4) circle (5pt);
\fill [color=black] (8,-4) circle (5pt);\fill [color=black] (9,-4) circle (5pt);
\fill [color=black] (10,-4) circle (5pt);\fill [color=black] (7,-5) circle (5pt);
\fill [color=black] (8,-5) circle (5pt);\fill [color=black] (9,-5) circle (5pt);
\fill [color=black] (10,-5) circle (5pt);\fill [color=black] (7,-6) circle (5pt);
\fill [color=black] (8,-6) circle (5pt);\fill [color=black] (9,-6) circle (5pt);
\fill [color=black] (10,-6) circle (5pt);
\draw (-6,6) circle (5pt);\fill (-5,6) circle (5pt);
\fill (-4,6) circle (5pt);\fill (-3,6) circle (5pt);
\fill (-2,6) circle (5pt);\fill (-1,6) circle (5pt);
\fill (-6,5) circle (5pt);\draw (-5,5) circle (5pt);
\fill (-4,5) circle (5pt);\fill (-3,5) circle (5pt);
\fill (-2,5) circle (5pt);\fill (-1,5) circle (5pt);
\fill (-6,4) circle (5pt);\fill (-5,4) circle (5pt);
\draw (-4,4) circle (5pt);\fill (-3,4) circle (5pt);
\fill (-2,4) circle (5pt);\fill (-1,4) circle (5pt);
\fill (-6,3) circle (5pt);\fill (-5,3) circle (5pt);
\fill (-4,3) circle (5pt);\draw (-3,3) circle (5pt);
\fill (-2,3) circle (5pt);\fill (-1,3) circle (5pt);
\fill (-6,2) circle (5pt);\fill (-5,2) circle (5pt);
\fill (-4,2) circle (5pt);\fill (-3,2) circle (5pt);
\fill (-1,2) circle (5pt);\draw (-2,2) circle (5pt);
\fill (-6,1) circle (5pt);\fill (-5,1) circle (5pt);
\fill (-4,1) circle (5pt);\fill (-3,1) circle (5pt);
\fill (-2,1) circle (5pt);\draw (-1,1) circle (5pt);
\fill(-1,0) circle (5pt);\fill (-2,0) circle (5pt);
\fill (-3,0) circle (5pt);\fill (-4,0) circle (5pt);
\fill (-5,0) circle (5pt);\fill (-6,0) circle (5pt);
\fill (-7,0) circle (5pt);\fill (-8,0) circle (5pt);
\fill (-9,0) circle (5pt);\fill (-10,0) circle (5pt);
\draw (-11,0) circle (5pt);\draw (-11,1) circle (5pt);
\draw (-11,2) circle (5pt);\draw (-11,3) circle (5pt);
\draw (-11,4) circle (5pt);\draw (-11,5) circle (5pt);
\draw (-11,6) circle (5pt);\draw (-11,7) circle (5pt);
\draw (-11,8) circle (5pt);\draw (-11,9) circle (5pt);
\draw (-11,10) circle (5pt);\draw (-11,11) circle (5pt);
\draw (-10,11) circle (5pt);\draw (-9,11) circle (5pt);
\draw (-8,11) circle (5pt);\draw (-7,11) circle (5pt);
\draw (-6,11) circle (5pt);\draw (-5,11) circle (5pt);
\draw (-4,11) circle (5pt);\draw (-3,11) circle (5pt);
\draw (-2,11) circle (5pt);\draw (-1,11) circle (5pt);
\fill [color=black] (-7,2) circle (5pt);
\fill (-6,-1) circle (5pt);\fill (-5,-1) circle (5pt);
\fill (-4,-1) circle (5pt);\fill (-3,-1) circle (5pt);
\fill (-2,-1) circle (5pt);\draw (-1,-1) circle (5pt);
\fill (-6,-2) circle (5pt);\fill (-5,-2) circle (5pt);
\fill(-4,-2) circle (5pt);\fill (-3,-2) circle (5pt);
\draw (-2,-2) circle (5pt);\fill (-1,-2) circle (5pt);
\fill (-6,-3) circle (5pt);\fill (-5,-3) circle (5pt);
\fill (-4,-3) circle (5pt);\draw (-3,-3) circle (5pt);
\fill (-2,-3) circle (5pt);\fill (-1,-3) circle (5pt);
\fill (-6,-4) circle (5pt);\fill (-5,-4) circle (5pt);
\draw (-4,-4) circle (5pt);\fill (-3,-4) circle (5pt);
\fill (-2,-4) circle (5pt);\fill (-1,-4) circle (5pt);
\fill (-6,-5) circle (5pt);\draw (-5,-5) circle (5pt);
\fill (-4,-5) circle (5pt);\fill (-3,-5) circle (5pt);
\fill (-2,-5) circle (5pt);\fill (-1,-5) circle (5pt);
\draw (-6,-6) circle (5pt);\fill (-5,-6) circle (5pt);
\fill (-4,-6) circle (5pt);\fill (-3,-6) circle (5pt);
\fill (-2,-6) circle (5pt);\fill (-1,-6) circle (5pt);
\fill (0,-1) circle (5pt);\draw (1,-1) circle (5pt);
\fill (2,-1) circle (5pt);\fill (3,-1) circle (5pt);
\fill (4,-1) circle (5pt);\fill (5,-1) circle (5pt);
\fill (6,-1) circle (5pt);\fill (0,-2) circle (5pt);
\fill (1,-2) circle (5pt);\draw (2,-2) circle (5pt);
\fill (3,-2) circle (5pt);\fill (4,-2) circle (5pt);
\fill (5,-2) circle (5pt);\fill (6,-2) circle (5pt);
\fill (0,-3) circle (5pt);\fill (1,-3) circle (5pt);
\fill (2,-3) circle (5pt);\draw (3,-3) circle (5pt);
\fill (4,-3) circle (5pt);\fill (5,-3) circle (5pt);
\fill (6,-3) circle (5pt);\fill (0,-4) circle (5pt);
\fill (1,-4) circle (5pt);\fill (2,-4) circle (5pt);
\fill (3,-4) circle (5pt);\draw (4,-4) circle (5pt);
\fill (5,-4) circle (5pt);\fill (6,-4) circle (5pt);
\fill (0,-5) circle (5pt);\fill (1,-5) circle (5pt);
\fill (2,-5) circle (5pt);\fill (3,-5) circle (5pt);
\fill (4,-5) circle (5pt);\draw (5,-5) circle (5pt);
\fill (6,-5) circle (5pt);\fill (0,-6) circle (5pt);
\fill (1,-6) circle (5pt);\fill (2,-6) circle (5pt);
\fill (3,-6) circle (5pt);\fill (4,-6) circle (5pt);
\fill (5,-6) circle (5pt);\draw (6,-6) circle (5pt);
\fill (0,-7) circle (5pt);\fill (0,-8) circle (5pt);
\fill (0,-9) circle (5pt);\fill (0,-10) circle (5pt);
\draw (-11,-1) circle (5pt);\draw (-11,-2) circle (5pt);
\draw (-11,-3) circle (5pt);\draw (-11,-4) circle (5pt);
\draw (-11,-5) circle (5pt);\draw (-11,-6) circle (5pt);
\draw (-11,-7) circle (5pt);\draw (-11,-8) circle (5pt);
\draw (-11,-9) circle (5pt);\draw (-11,-10) circle (5pt);
\draw (-11,-11) circle (5pt);\draw (-10,-11) circle (5pt);
\draw (-9,-11) circle (5pt);\draw (-8,-11) circle (5pt);
\draw (-7,-11) circle (5pt);\draw (-6,-11) circle (5pt);
\draw (-5,-11) circle (5pt);\draw (-4,-11) circle (5pt);
\draw (-3,-11) circle (5pt);\draw (-2,-11) circle (5pt);
\draw (-1,-11) circle (5pt);\draw (0,-11) circle (5pt);
\draw (1,-11) circle (5pt);\draw (2,-11) circle (5pt);
\draw (3,-11) circle (5pt);\draw (4,-11) circle (5pt);
\draw (5,-11) circle (5pt);\draw (6,-11) circle (5pt);
\draw (7,-11) circle (5pt);\draw (8,-11) circle (5pt);
\draw (9,-11) circle (5pt);\draw (10,-11) circle (5pt);
\draw (11,-11) circle (5pt);\draw (11,-10) circle (5pt);
\draw (11,-9) circle (5pt);\draw (11,-8) circle (5pt);
\draw (11,-7) circle (5pt);\draw (11,-6) circle (5pt);
\draw (11,-5) circle (5pt);\draw (11,-4) circle (5pt);
\draw (11,-3) circle (5pt);\draw (11,-2) circle (5pt);
\draw (11,-1) circle (5pt);
\end{tikzpicture}
\caption{The lattice points corresponding to $(k,l)$-colorings of $\pm K_2$ with $k=6$ and $l=4$.}
\label{polytopalS}
\end{figure}
Thus, by \eqref{iopreciprocity} and \eqref{phirecsignedeq},
\begin{align*}
  (-1)^{ |V| } c_\Sigma (-2k-1, -2l)
  &= (-1)^{ |V| } L_{ \Box^\circ, \HH_\Sigma }^\circ (-k) + (-1)^{ |V| } \Phi(-k-1,-l) \\
  &= L_{ \Box, \HH_\Sigma } (k) + \left| (k+l) \Box \cap \Z^V \right| - \left| k \Box \cap \Z^V \right| .
\end{align*}
What we are counting on the right-hand side are the integer lattice points in the cube $[-k-l,k+l]^V$, with multiplicity equal to $1$ if outside the cube $[-k,k]^V$, otherwise with multiplicity equal to the number of closed regions of $\HH_\Sigma$ the points lies in. As shown in \cite{zaslavskysignedcoloring} (and again used in \cite{iop}), the latter can be interpreted as the number of compatible acyclic orientations of~$\Sigma$.
\end{proof}



\section{Open Questions}\label{sec:openprob}

We finish with two venues for future research.
First, one can associate several matroids to a signed graph, most prominently Zaslavsky's \emph{frame matriod} and
\emph{(extended) lift matroid} \cite{zaslavskysignedgraphs,zaslavskydynamicsurvey}. It is a natural question to ask about possible connections between the Tutte polynomials of these matroids and the bivariate chromatic polynomials.
Second, \emph{gain} and \emph{biased graphs} are natural generalizations of signed graphs \cite{zaslavskydynamicsurvey}, and
so another natural question concerns possible extensions of our work to these more general constructs.

\bibliographystyle{amsplain}

\def\cprime{$'$} \def\cprime{$'$}
\providecommand{\bysame}{\leavevmode\hbox to3em{\hrulefill}\thinspace}
\providecommand{\MR}{\relax\ifhmode\unskip\space\fi MR }
\providecommand{\MRhref}[2]{%
  \href{http://www.ams.org/mathscinet-getitem?mr=#1}{#2}
}
\providecommand{\href}[2]{#2}

\setlength{\parskip}{0cm} 

\end{document}